\documentclass{article}
\pdfoutput=1
\usepackage{amssymb, amsmath ,graphicx,
paralist,cancel,caption,subcaption, amsthm}
\usepackage[numbers,sort]{natbib}
\usepackage[titletoc,title]{appendix}
\usepackage{algorithm, algpseudocode, fullpage}
\numberwithin{algorithm}{section}

\usepackage{tikz}
\usepackage{pgfplots}
\pgfplotsset{compat=1.3} 
\usetikzlibrary{calc}
\usetikzlibrary{patterns}
\usepackage{listings, color}
\usepackage{cleveref}

\newcommand{\B}{{\bf B}}
\newcommand{\U}{{\bf U}}
\newcommand{\bPhi}{{\bf \Phi}}
\newcommand{\bPsi}{{\bf \Psi}}

\newcommand{\E}{{\bf E}}
\newcommand{\prop}{{\mathcal P}}
\newcommand{\SO}{{\mathcal O}}
\newcommand{\R}{{\mathbb R}}
\newcommand{\ZZ}{{\mathbb Z}}

\newcommand{\SL}{{\mathcal L}}
\newcommand{\J}{{\bf J}}

\renewcommand{\k}{{\bf k}}

\newcommand{\x}{{\bf x}}
\newcommand{\z}{{\bf z}}

\newcommand{\f}{{\bf f}}

\newcommand{\g}{{\bf g}}
\newcommand{\h}{{\bf h}}
\newcommand{\y}{{\bf y}}
\newcommand{\K}{{\bf K}}
\renewcommand{\v}{{\bf v}}

\newcommand{\samp}{{\mathcal S}}
\newcommand{\interp}{{\mathcal I}}

\begin{document}
\title{An adaptive local discrete convolution method for the numerical solution of Maxwell's equations}
\author
{
	Boris Lo\thanks{Graduate Group in Applied Science and Technology, University of California, Berkeley, Berkeley, CA 94720.  Email \texttt{bt.lo@berkeley.edu}}\and
	Phillip Colella\thanks{Lawrence Berkeley National Laboratory, Berkeley, CA 94720. Email: \texttt{pcolella@lbl.gov}}
}

\maketitle
\begin{abstract}
We present a numerical method for solving the free-space Maxwell's equations in three dimensions using compact convolution kernels on a rectangular grid. We first rewrite Maxwell's Equations as a system of wave equations with auxiliary variables and discretize its solution from the method of spherical means. The algorithm has been extended to be used on a locally-refined nested hierarchy of rectangular grids. 
\end{abstract}

\section{Introduction}\label{sec:intro}

We want to solve the free-space 3D Maxwell's equations
\begin{gather}
\frac{\partial\E}{\partial t} = c\nabla\times\B - 4\pi\J,\\
\frac{\partial\B}{\partial t} = -c\nabla\times\E, \\
\nabla\cdot\E = 4\pi\rho,\\
\nabla\cdot\B = 0.
\end{gather}

In our previous work \cite{lo}, we considered Maxwell's equations in Fourier space, derived a real space propagator for the system, and discretized the exact solution from Duhamel's formula. This propagator includes Helmholtz decomposition operators. The Helmholtz decomposition operators require global Poisson solves at every time step which offsets the computational advantages of the local convolution kernel parts of the propagator. \\

In the present work, we get around this difficulty by applying a similar technique to an auxiliary system of equations instead of directly to Maxwell's equations. This auxiliary system is a system of wave equations for $\E,\B$ combined with constraints which, if satisfied initially, are satisfied for all time, such that the solutions of the auxiliary system are solutions to Maxwell's equations. We then apply Kirchhoff's formula to this system and discretize the resulting convolution equations. The convolution kernels from this propagator are the same as the local kernels for the transverse Maxwell's equations' propagator in \cite{lo} and thus the same discretization techniques and domain decomposition can be applied. The locality of the convolution kernels allows us to naturally incorporate adaptive mesh refinement (AMR), where the domain is divided up in nested hierarchy of rectangular grids at each refinement level. \\

In Section \ref{sec:outline} we introduce the auxiliary system and show the analytic solution for Maxwell's equations in terms of a propagator with specified charges and currents. In Section \ref{sec:discretization}, we describe the discretization process briefly, and discuss in detail the local discrete convolution method (LDCM) Maxwell solver for a single level and its extension to multiple levels. In Section \ref{sec:experiments} we present a number of numerical tests that show an implementation of our algorithm.  Finally, in Section \ref{sec:conclusion} we make some concluding remarks.

\section{Problem Statement and Derivation of Propagators}\label{sec:outline}


\subsection{Maxwell's Equations}\label{sec:maxwellProp}
Introducing $\bPhi \equiv \nabla \times \B$ and $\bPsi \equiv \nabla \times \E$, we rewrite Maxwell's Equations, with $\rho,\J$ specified, as the following auxiliary system of wave equations
\begin{gather}
\frac{\partial \E}{\partial t} = c \bPhi - 4\pi \J, \label{eq:Eeq}\\
\frac{\partial \bPhi}{\partial t} = c \nabla^2 \E - 4\pi c\nabla \rho, \label{eq:Phieq}\\
\frac{\partial \B}{\partial t} = -c \bPsi, \label{eq:Beq}\\
\frac{\partial \bPsi}{\partial t}  = -c \nabla^2 \B - 4\pi \nabla \times \J. \label{eq:Psieq}
\end{gather}
If the initial conditions satisfy 
\begin{gather}
\bPsi = \nabla \times \E, \label{eq:curlEconstraint} \\
\bPhi = \nabla \times \B, \label{eq:curlBconstraint} \\
\nabla\cdot\E = 4\pi\rho, \label{eq:divEconstraint}\\
\nabla\cdot\B=0, \label{eq:divBconstraint}
\end{gather}
then the auxiliary system is equivalent to the original Maxwell system. To show this, consider the four error quantities associated with the initial value constraints at $t=0$
\begin{gather}
\K_B = \bPhi - \nabla\times\B, \label{eq:erroreq1}\\
\K_E = \bPsi - \nabla\times\E, \label{eq:erroreq2}\\
D_B = \nabla\cdot\B, \label{eq:erroreq3}\\
D_E = \nabla\cdot\E - 4\pi\rho. \label{eq:erroreq4}
\end{gather}
Using the auxiliary system \eqref{eq:Eeq}-\eqref{eq:Psieq}, the four evolution equations associated with these quantities are given by
\begin{gather}
\frac{\partial\K_B}{\partial t} = c\nabla\times\K_E + c\nabla D_E, \label{eq:errorevo1}\\
\frac{\partial\K_E}{\partial t} = -c\nabla\times\K_B - c\nabla D_B, \\
\frac{\partial D_B}{\partial t} = -c\nabla\cdot\K_E, \\
\frac{\partial D_E}{\partial t} = c\nabla\cdot\K_B. \label{eq:errorevo4}
\end{gather}
It is clear that if $\K_B,\K_E,D_B,D_E$ vanish at $t=0$, then they remain zero for all time after. In particular, the symbol of the linear operator associated with these eight evolution equations has the eigenvalues $\pm ic|\k|$ each with a multiplicity of four. Since errors propagate away with the same wave speed, any error will not accumulate at a fixed location and be a potential source of numerical instability. The initial value problem \eqref{eq:Eeq}-\eqref{eq:Psieq} is well-posed even if the initial-value constraints \eqref{eq:erroreq1}-\eqref{eq:erroreq4} are not satisfied. The constraints are required only so that the solution is equivalent to the solution to Maxwell's Equations. Since the two systems are equivalent, the solutions for $\E,\B$ obtained from the auxiliary system will also be the solution to the original Maxwell system. 

The solutions to \eqref{eq:Eeq}-\eqref{eq:Psieq} are given by Kirchhoff's Formula using the method of spherical means \cite[p.231]{whitham}. Defining the kernels $G^{\Delta t}$ and $H^{\Delta t}$ as
\begin{gather}
G^{\Delta t}(\z) \equiv \frac{\delta(|\z|-c\Delta t)}{4\pi c\Delta t}, \label{eq:deltaG}\\
H^{\Delta t}(\z) \equiv \frac{1}{c}\frac{\partial}{\partial s}\left(\frac{\delta(|\z|-cs)}{4\pi cs}\right)\bigg|_{s=\Delta t}. \label{eq:deltaH}
\end{gather}
$G^{\Delta t}$ is a spherical delta distribution with radius $c\Delta t$. The action of the propagator on an arbitrary state vector $\h(\x) \equiv[\f(\x),\, \g(\x)]^T$ with $\f,\g\in\R^3$ is given by
\begin{gather}\label{eq:propagator}
\prop^{\Delta t}[\h] = \left[\begin{array}{c} H^{\Delta t} * \f + G^{\Delta t} * \g\\ G^{\Delta t} * \nabla^2\f+  H^{\Delta t} * \g\end{array} \right],
\end{gather}
where the scalar convolution kernel with vector quantity is defined as convolution with each component and that convolutions are defined spatially as
\begin{gather}
(K*f)(\x) \equiv \int_{\R^3} K(\y) f(\x-\y)\,d\y.
\end{gather}
In particular, the solution to \eqref{eq:Eeq}-\eqref{eq:Phieq} is then given by
\begin{gather}\label{eq:Esol}
\left(\begin{array}{c}\E(\x,t+\Delta t) \\ \bPhi(\x,t+\Delta t) \end{array}\right) = \prop^{\Delta t}\left[\left(\begin{array}{c}\E(\x,t) \\ \bPhi(\x,t) \end{array}\right)\right] - 4\pi \int_{t}^{t+\Delta t} \prop^{t+\Delta t-s}\left[\left(\begin{array}{c}\J(\x,s) \\ c\nabla\rho(\x,s) \end{array}\right)\right] ds.
\end{gather} 
The propagator for \eqref{eq:Beq}-\eqref{eq:Psieq} is the same as that for \eqref{eq:Eeq}-\eqref{eq:Phieq}, with the substitution $\Delta t\rightarrow-\Delta t$. Thus, the solution is given by
\begin{gather}\label{eq:Bsol}
\left(\begin{array}{c}\B(\x,t+\Delta t) \\ \bPsi(\x,t+\Delta t) \end{array}\right) = \prop^{-\Delta t}\left[\left(\begin{array}{c}\B(\x,t) \\ \bPsi(\x,t) \end{array}\right)\right] - 4\pi \int_{t}^{t+\Delta t} \prop^{-(t+\Delta t-s)}\left[\left(\begin{array}{c}0 \\ \nabla\times\J(\x,s) \end{array}\right)\right] ds.
\end{gather}
It can be seen from the Fourier transforms of the convolution kernels that
\begin{gather}
G^{-\Delta t}*f = -G^{\Delta t}*f, \\
H^{-\Delta t}*f = H^{\Delta t}*f. 
\end{gather}
In addition
\begin{gather}
H^{\Delta t}*f = \frac{1}{ct}G^{\Delta t}*f - \sum_{i=1}^3 G_i^{\Delta t}*\frac{\partial f}{\partial z_i}, \label{eq:Hlinear} \\
G_i^{\Delta t}(\z)  = \frac{z_i\delta(|\z|-c\Delta t)}{4\pi c \Delta t}.
\end{gather}
With these, we have fully specified the solutions, \eqref{eq:Esol} and \eqref{eq:Bsol}, in terms of convolution with weighted spherical delta distributions. We note that it can be shown directly that $\bPsi(\x,t+\Delta t) = \nabla\times\E(\x,t+\Delta t)$ and $\bPhi(\x,t+\Delta t) = \nabla\times\B(\x,t+\Delta t)$ given the constraints are satisfied at $t$. When $\rho,\J$ are not a specified but functions of field variables, instead of using Kirchhoff's formula and a quadrature scheme one can use Lawson's method \cite{lawson} for time integration.

\section{Discretization Approach}\label{sec:discretization}

\subsection{Single Level Algorithm}
We consider a rectangular domain discretized with a Cartesian grid with grid spacing $h$ with open boundary conditions. The convolutions in \eqref{eq:Esol}-\eqref{eq:Bsol} are approximated with discrete convolutions on the grid. This requires a discretized representation of the convolution kernels, $G^{\Delta t,h}\approx G^{\Delta t}(\z),H^{\Delta t,h}\approx H^{\Delta t}(\z)$, on the grid. $H^{\Delta t,h}$ is obtained by \eqref{eq:Hlinear}, so that the problem reduces to only creating discrete representations of (weighted) spherical delta distributions. We refer the reader to \cite{lo} for a detailed treatment of the discretization of the convolution kernels. The resulting discrete convolution kernels have compact support just like their continuous counterparts. Thus, the discrete convolutions can be computed exactly using Hockney's method \cite{hockney}. \\

The overall time-stepping algorithm is given in Algorithm \ref{alg:singlelevel_alg}. This defines the discrete evolution for $\E,\B$, since $\bPhi,\bPsi$ are computed at the beginning of every time step. The source term integrals are discretized using a closed Newton-Cotes quadrature scheme with step size $\Delta s = \Delta t/(M-1)$ where $M$ is the number of quadrature points.  We choose a fixed step size quadrature because $\prop^{t_1}[\prop^{t_2}[\U]]=\prop^{t_1+t_2}[\U]$ and therefore we only need to create one propagator with step size $\Delta s$ during initial setup. \\

Even though the divergence constraints are preserved by the continuous time evolution, deviations from \eqref{eq:divEconstraint}-\eqref{eq:divBconstraint} may be generated by discretization error.  To help remedy this, we apply local filters \cite{marder} of the form
\begin{align}\label{eq:filter}
\E &:= \E + \eta(\SL\E - 4\pi\nabla\rho), \\
\B &:= \B + \eta\SL\B, \\
\SL_{ij} &= \partial_{x_i}\partial_{x_j},
\end{align}
where $\eta\sim \SO(h^2)$ is a constant and $\SL$ is a matrix valued operator with the diagonal terms discretized with centered-difference approximations to the second derivative while the off-diagonal terms are products of centered-difference approximations to the first derivatives. This filtering step corresponds to applying an explicit diffusion step to the error in the longitudinal fields.  Note that we do not have to do this for the curl constraints \eqref{eq:curlEconstraint}-\eqref{eq:curlBconstraint}, since $\bPhi,\bPsi$ are re-initialized at the beginning of each time step.

\begin{algorithm}[h!]
  \caption{Single level LDCM for Maxwell's Equations}
  \label{alg:singlelevel_alg}
  \begin{algorithmic}
  \State Initialize Newton-Cotes quadrature weights $\{w_m\}_{m=0}^M$\
  \State \texttt{/* Create the convolution kernels with quadrature step size $\Delta s$ and spacing $h$ */}
  \State Compute $G^{\Delta s,h}$, and $H^{\Delta s,h}$ 
    \State \texttt{/* Begin time-stepping loop */}
      \For{$n = 1,2,\dots$}
        \State \texttt{/* Initialize the fields for this time step */}
        \State \texttt{/* Let $U^{(n),h}\approx U(n\Delta t,\x)$ */}
        \State $\E^{(n),h} \gets \E^{(n-1),h}, \B^{(n),h} \gets \B^{(n-1),h}, \bPhi^{(n),h} \gets\nabla\times\E^{(n),h}, \bPsi^{(n),h}\gets\nabla\times\B^{(n),h}$
      \State \texttt{/* Begin quadrature loop */}
        \For{$m = 1,2,\dots,M$}
            \State \texttt{/* Add in source terms evaluated at $t=(n-1)\Delta t+(m-1)\Delta s$ */}
            \State $\E^{(n),h} \gets \E^{(n),h} - w_m 4\pi\J^h$
            \State $\bPhi^{(n),h} \gets \bPhi^{(n),h} - w_m 4\pi c\nabla\rho^h$
            \State $\bPsi^{(n),h} \gets \bPsi^{(n),h} - w_m 4\pi\nabla\times\J^h$
            \State \texttt{/* Apply propagator to the fields except final quadrature point */}
            \If{$m<M$}
            \State $\left[\begin{array}{c}\E^{(n),h} \\ \bPhi^{(n),h}\end{array}\right] \gets \left[\begin{array}{c}H^{\Delta s,h} * \E^{(n),h} + G^{\Delta s,h}*\bPhi^{(n),h} \\ (G^{\Delta s,h}*\nabla^2) * \E^{(n),h} + H^{\Delta s,h}*\bPhi^{(n),h}\end{array}\right]$
            \State $\left[\begin{array}{c}\B^{(n),h} \\ \bPsi^{(n),h}\end{array}\right] \gets \left[\begin{array}{c}H^{\Delta s,h} * \B^{(n),h} - G^{\Delta s,h}*\bPsi^{(n),h} \\ -(G^{\Delta s,h}*\nabla^2) * \B^{(n),h} + H^{\Delta s,h}*\bPsi^{(n),h}\end{array}\right]$
            \EndIf
        \EndFor
        \State \texttt{/* Enforcing Constraints */}
        \State $\E^{(n),h}\leftarrow\E^{(n),h}+\eta(\SL\E^{(n),h} - 4\pi\nabla\rho^{h})$ 
        \State $\B^{(n),h}\leftarrow\B^{(n),h}+\eta\SL\B^{(n),h}$ 
      \EndFor 
  \end{algorithmic}
 \end{algorithm}

\subsection{Domain Decomposition}
Since the discretized version of the propagator involves only local operators, we can use standard domain decomposition to parallelize this algorithm.  Consider a single level domain, $\Omega_h$, partitioned into rectangular patches. For each patch:  
\begin{enumerate}
\item at the beginning of each quadrature step, copy field values in ghost region from neighboring processors,
\item apply propagator to update local field values, invalidating values in ghost region.
\end{enumerate}
The minimum width of the ghost region is determined by the size of the quadrature, $\Delta s$, and the order of the method because the size of the support of the spherical delta distributions is dependent on how far in time the fields are to be advanced. \\

For a point, $\x_k$, near the boundary, when applying the discrete convolutions we replace the field values outside the computational domain with the current field value at $\x_k$. This approximation leads to waves reflecting back into the computational domain. We mitigate this error with mesh refinement, by placing the boundary of the computational domain far away from the sources. The amplitude of the waves reaching the boundary will be weaker and the reflected error waves will also be smaller. 

\subsection{Multilevel Algorithm}
\begin{figure}
\centering
\begin{tikzpicture}[scale=0.3] 
\draw[step=1cm] (2,2) grid (10,10);
\draw[step=0.5cm, fill=black!20!white] (-5,-3) grid (-1,1) rectangle (-5, -3);
\draw[step=0.5cm, fill=white] (-4,-2) grid (-2,0) rectangle (-4,-2);
\draw[ultra thick] (-4,-2) rectangle (-2,0);
\draw[ultra thick] (5,5) rectangle (7,7);
\draw[dashed] (-4,-2) -- (5,5);
\draw[dashed] (-4,0) -- (5,7);
\draw[dashed] (-2,-2) -- (7,5);
\draw[dashed] (-2,0) -- (7,7);
\draw (2,-2.5) node {\large $\Omega_1\cup\Omega_{1,g}$};
\draw (10.9,2.5) node {\large$\Omega_0$};
\end{tikzpicture}
\caption{Example schematic of a two-level nested domains with factor of 2 refinement. The unshaded region is $\Omega_1$ and the shaded region is the ghost region $\Omega_{1,g}$.}
\label{fig:multilevel_cartoon}
\end{figure}
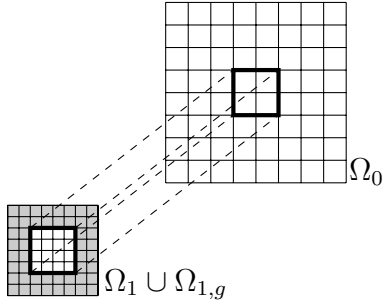
Consider now a hierarchy of nested rectangular grids, $\Omega_j,\;j=0,\ldots,J-1$, where the grid spacing for $\Omega^j$ is $\frac{h}{r^j}$ for some refinement factor, $r\in\ZZ^+$, with $\Omega_j\cup\Omega_{j+1}=\Omega_{j+1}$, $j=0\ldots,J-2$.  We introduce sampling and interpolation operator, $\samp$ and $\interp$ respectively, to communicate field values with the next immediate lower and upper levels. Similar to the ghost regions for each patch in parallelizing the single level algorithm, we define a ghost region for each level, $\Omega_{j,g}$, where the width of the ghost region is determined by how far in the time fields are to be advanced. At the beginning of each quadrature step, except on the first level, for all nodes in $\Omega_{j,g}$ we interpolate $\E,\B$ from level $j-1$.  After interpolating, except on the finest level, we replace $\E,\B$ at level $j$ with field values from level $j+1$ on the nodes that are in $\Omega_j\cap\Omega_{j+1}$. A sample schematic of two levels with $r=2$ is shown in Figure \ref{fig:multilevel_cartoon}.  After interpolating and sampling, each level is evolved independently with the propagator. \\

Let $f_j^{(n)}$ denote discretized $f$ on level $j$ and at time $t_n=n\Delta t$, the multilevel algorithm is outlined in Algorithm \ref{alg:multilevel_alg}. Since \eqref{eq:Eeq}-\eqref{eq:Psieq} is a system of linear differential equations, we can use linear superposition to generate the overall solution to the problem in this multilevel setup; the solution is given by a composite where it takes the finest level values for any subdomain. For example, in the two level case, let $\U=(\E,\B,\bPsi,\bPhi)^T$, then the solution is given by
\begin{align}
\U^{(n)} = \left\{\begin{array}{l l} \U^{(n)}_1 ,& \text{on }\Omega_1 \\
\U^{(n)}_0 ,& \text{on }\Omega_0\setminus\Omega_1
\end{array} \right. .
\end{align}

\begin{algorithm}[h!]
  \caption{Multilevel LDCM for Maxwell's Equations}
  \label{alg:multilevel_alg}
  \begin{algorithmic}
  \State Initialize Newton-Cotes quadrature weights $\{w_m\}_{m=0}^M$
  \ForAll {levels $\Omega_j,\;j=0,\ldots,J-1$}
    \State Initialize $\U_j^{(0)}$
    \State Compute $G^{\Delta s,h/r^j}$, and $H^{\Delta s,h/r^j}$ 
  \EndFor
    \State \texttt{/* Begin time-stepping loop */}
      \For{$n = 1,2,\dots$}
        \ForAll {levels $\Omega_j,\;j=0,\ldots,J-1$}
        \State \texttt{/* Initialize the fields for this time step */}
        \State $\U_j^{(n)} \gets \U_j^{(n-1)}$    
        \For{quadrature step $s$}
        \State \texttt{/* Apply sampling operator except for level 0 */}
        \State $\U_{j-1}^{(n)}\leftarrow\samp[\U_{j}^{(n)}]$ on $\Omega_{j}$
        \State \texttt{/* Apply interpolation operator except for level $J-1$ */}
        \State $\U_{j+1}^{(n)}\leftarrow\interp[\U_{j}^{(n)}]$ on $\Omega_{j+1,g}$
        \State Apply single level operations (add in source term and apply propagator)
        \EndFor
      \EndFor 
      \State Sample and interpolate $\E,\B$ so that $\SL$ can be applied on the refinement levels
      \State Enforce the constraints independently for each level 
  \EndFor
  \end{algorithmic}
 \end{algorithm}

Since we interpolate once every quadrature step, the width of $\Omega_{i,g}$ for level $i$ has the same width as the ghost region required for domain decomposition.

\subsubsection{Interpolation}
We use high order B-splines (see Appendix \ref{sec:appendix}) to interpolate the fields between levels similar to the ones used to regularize the delta distributions in the propagator. However, the choice of interpolant is more restrictive than the one used to regularize the delta distribution. The convergence of spherical quadrature when regularizing the delta distribution depends on the smoothness of the integrand \cite{atkinson}. However, we are interested in the regularized delta distribution as a discrete convolution kernel with some discretized function $f$. Numerically, the spherical quadrature and discrete convolution commutes and therefore we relied on the smoothness of $f$ for the convergence of the spherical quadrature. This allows us to use a $C^0$ high order B-spline as a regularizer with the advantage that it has minimal support. \\

In this method, $f$ is a field component or a component of the source terms. Since the field components must be sufficiently smooth for the spherical quadrature and the accuracy of the high order finite difference operators applied to the field components also depend on smoothness, these translate into a smoothness requirement for the interpolants. For a $q$-th order method, we would need the error from the spherical quadrature to be at least $\mathcal{O}(h^q)$ which requires $f\in C^q$.  Therefore the interpolant must also be at least $q$-th order accurate and $C^q$.

\subsubsection{Regridding}
 For an adaptive version of this method, instead of a fixed hierarchy of rectangular grids, we regrid at the beginning of any time step as needed.  Suppose we wish to regrid level $j$, $j>0$, let $\Omega_j=\Omega_{j,discard}\cup\Omega_{j,keep}$ before regridding and $\Omega_j=\Omega_{j,keep}\cup\Omega_{j,new}$ after regridding. First sample down on $\Omega_{j,discard}$, then interpolate on $\Omega_{j,new}$ using the same sampling and interpolating operators. The regridding algorithm is outlined in Algorithm \ref{alg:regrid_alg}.

\begin{algorithm}[h!]
  \caption{Regridding Algorithm}
  \label{alg:regrid_alg}
  \begin{algorithmic}
  \For{levels $\Omega_j,\;j=1,\ldots,J-1$}
    \If{regrid}
      \State \texttt{/* Sample down starting from topmost level */}
      \For{$k=J-1,\ldots,j$}
        \State $\U_{k-1}^{(n)}\leftarrow\samp[\U_{k}^{(n)}]$ on $\Omega_{k}\cap \Omega_{j,discard}$
        \State \texttt{/* Discard part of domain that has been sampled from */}
        \State $\Omega_k \leftarrow \Omega_k \setminus (\Omega_{k}\cap \Omega_{j,discard})$
      \EndFor
      \State $\Omega_j \leftarrow \Omega_j \cup \Omega_{j,new}$
      \State \texttt{/* Interpolate from level $j-1$ */}
        \State $\U_{j}^{(n)}\leftarrow\interp[\U_{j-1}^{(n)}]$ on $\Omega_{j,new}$
      \State Enforce the constraints
    \EndIf
  \EndFor
  \end{algorithmic}
 \end{algorithm}

\section{Numerical Results}\label{sec:experiments}

We implemented a fourth order version of our Maxwell solver with $c=1$; the one step error for the solver is $\SO(h^5)$ but after some number of time steps the total error will be $\SO(h^{q-1})$ for a method that has a one step error of $\SO(h^{q})$ and $\Delta t=\SO(h)$. We used sixth-order centered-differences for the spatial derivatives, the fifth-order 3/8 Simpson's Rule for the source integration, $W_{6,0}$ for the discrete delta distribution, and $W_{6,6}$ for the interpolation operator. The discrete convolutions are performed via Hockney's method extending the domain equal to the support of the discrete convolution kernels and using the FFTW library \cite{fftw}. The domain at the coarsest level is a unit cube and each level is divided into $33^3$ node patches with factor of four refinement; every level has the same number of nodes, $N$.  The filter parameter at level $j$ is $\eta_j=\frac{45}{544}h_j^2$. For each test, $\Delta t$ is the same across refinement levels. 
  
\subsection{Translating Spherical Charge Distribution}
\begin{figure*}
        \centering
        \begin{subfigure}[t]{0.475\textwidth}
            \centering
            \includegraphics[width=\textwidth]{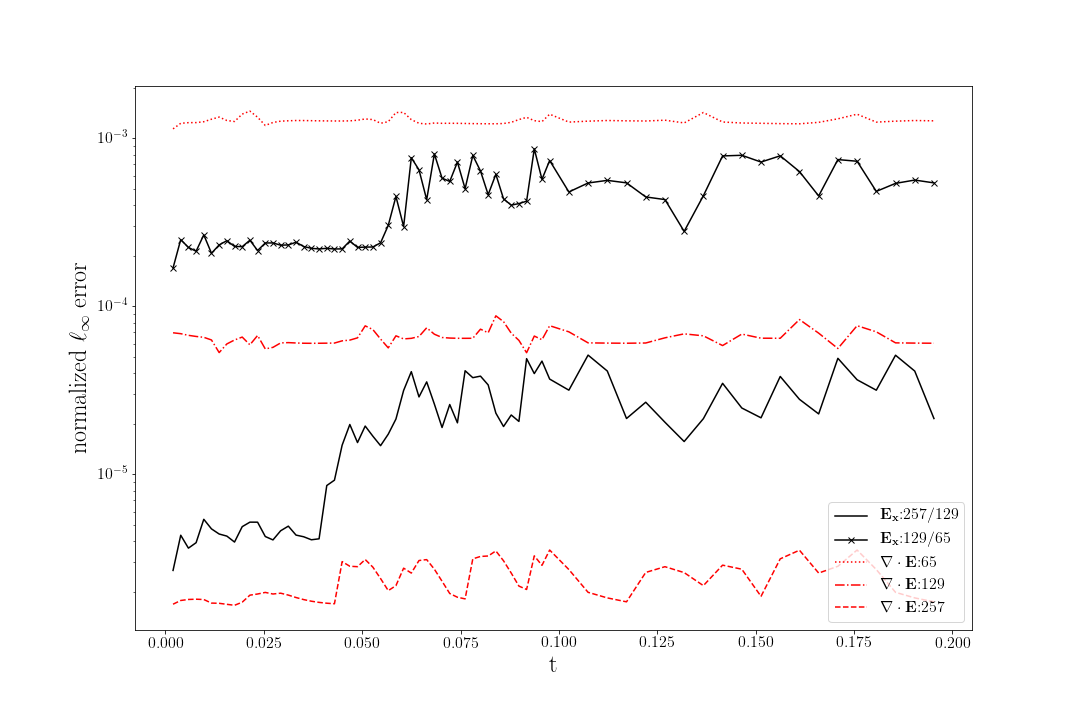}  
        \end{subfigure}
        \hfill
        \begin{subfigure}[t]{0.475\textwidth}  
            \centering 
            \includegraphics[width=\textwidth]{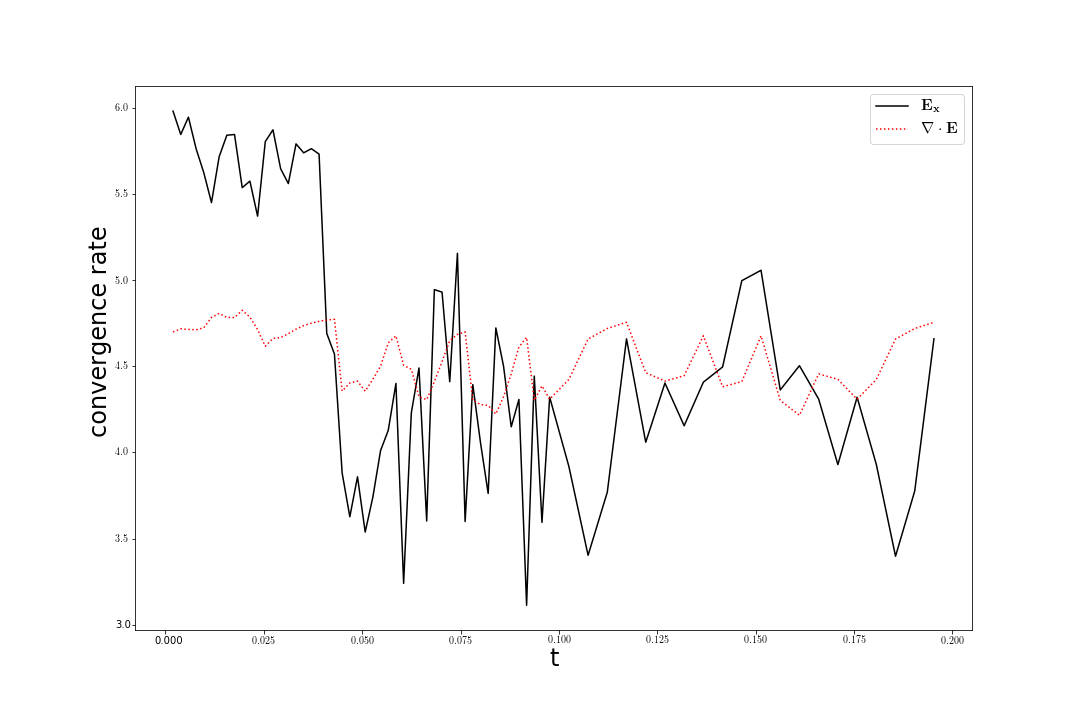}
        \end{subfigure}
        \caption[]
        {$\ell_{\infty}$ error values and convergence results for $E_x$ and $\nabla\cdot\E-4\pi\rho$ for the translating spherical charge distribution problem as a function of time in $\Omega_2$. On the left are the normalized $\ell_{\infty}$ errors for $E_x$ and $\nabla\cdot\E-4\pi\rho$. The errors for $E_x$ are obtained from the difference of sampled field values from the $N=257$ with $N=129$ and also from sampled $N=129$ with $N=65$ test case. The $E_x$ error is normalized by the max norm of the electrostatic solution ($\approx 0.0694795$) and $\nabla\cdot\E-4\pi\rho$ error is normalized by $\max_{\x} 4\pi\rho\approx 30.6796$. On the right are the associated convergence rates.} 
        \label{fig:mc_real_pics}
\end{figure*}

For the first numerical test, we used a $C^6$ spherical support charge distribution with a spatially constant $\v(t)$.
\begin{align}
\rho(\x,t) &= \left\{\begin{array}{rl} a(r(t)-r(t)^2)^6,& r < 1 \\ 0,& r \geq 1 \end{array}\right. , \,\, r = \frac{1}{R_0}||\x-\x_0||,\\
\J(\x,t) &= \v(t)\rho(\x,t), \\
\v(t) &= \nu d \pi\frac{35}{16}\sin^7(2\pi\nu t) \hat\v.
\end{align}
The electrostatic solution is given by
\begin{align}
\E(\x) &= 4\pi R_0 a\hat r\left\{\begin{array}{rl} \frac{r^7}{9} - \frac{3r^8}{5} + \frac{15r^9}{11} - \frac{5r^{10}}{3} + \frac{15r^{11}}{13} - \frac{3r^{12}}{7} + \frac{r^{13}}{15},& r < 1 \\ \frac{1}{45045r^2},& r \geq 1 \end{array}\right. ,\label{eq:spherical_es}\\
\B(\x) &= 0 .
\end{align}
$\hat r$ is with respect to $\x_0$ and we use this as the initial condition for this test problem. We perform this test on fixed grids with two refinement levels: $\Omega_1=\left[\frac{3}{8},\frac{5}{8}\right]^3$ and $\Omega_2=\left[\frac{15}{32},\frac{17}{32}\right]^3$, with parameters: $a=10^{4}, d=\frac{1}{256},\nu = \frac{1024}{80}$, $R_0=\frac{1}{72}$, $\x_0 = \left(\frac{127}{256}, \frac{127}{256}, \frac{127}{256}\right)$, $\hat\v=\left(\cos\frac{\sqrt{3}}{3}\cos\frac{\sqrt{2}}{3}, \sin\frac{\sqrt{3}}{3}\cos\frac{\sqrt{2}}{3}, \sin\frac{\sqrt{2}}{3}\right)$, $N=(65,129,257)$ with $\Delta t = \left(\frac{1}{1024}, \frac{1}{2048}, \frac{1}{4096}\right)$ respectively, this corresponds to CFL = 1 at the finest level, out to $t_{final}=\frac{200}{1024}$. Fig. \ref{fig:mc_real_pics} shows the $E_x$ Richardson convergence rate estimate and the associated $\ell_{\infty}$ error as well as the absolute convergence rate and associated $\ell_{\infty}$ errors for $\nabla\cdot\E-4\pi\rho$ on the three grids in $\Omega_2$ as a function of time step and as expected our solution shows fourth order convergence.

\subsubsection{Electrostatic Test}
\begin{figure*}
        \centering
        \begin{subfigure}[t]{0.475\textwidth}
            \centering
            \includegraphics[width=\textwidth]{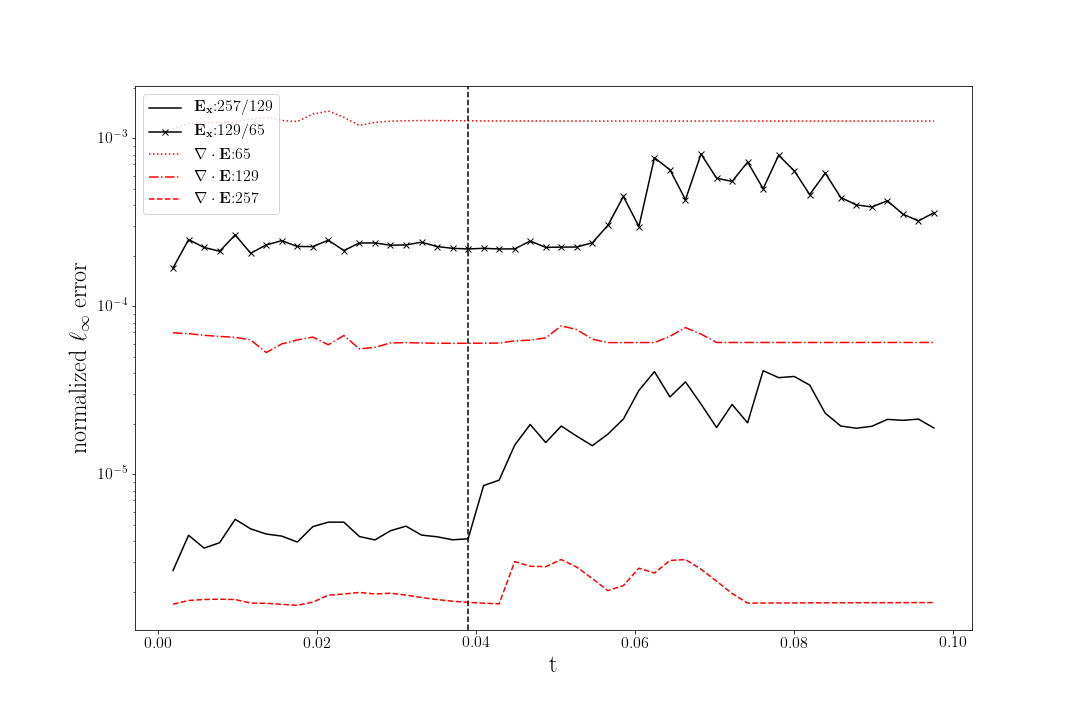}  
        \end{subfigure}
        \hfill
        \begin{subfigure}[t]{0.475\textwidth}  
            \centering 
            \includegraphics[width=\textwidth]{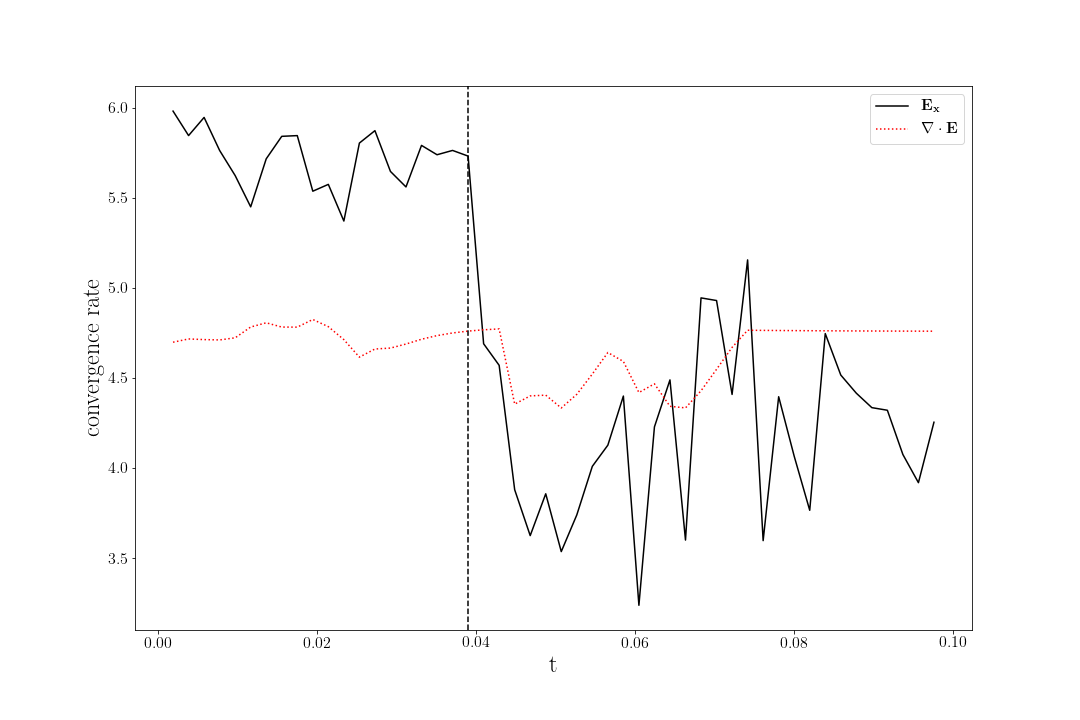}   
        \end{subfigure}
        \caption[]
        {$\ell_{\infty}$ error values and convergence results for $E_x$ and $\nabla\cdot\E-4\pi\rho$ for the stopped spherical charge distribution problem as a function of time in $\Omega_2$. On the left are the normalized $\ell_{\infty}$ errors for $E_x$ and $\nabla\cdot\E-4\pi\rho$. The errors for $E_x$ are obtained from the difference of sampled field values from the $N=257$ with $N=129$ and also from sampled $N=129$ with $N=65$ test case. The $E_x$ error is normalized by the max norm of the electrostatic solution ($\approx 0.0694795$) and $\nabla\cdot\E-4\pi\rho$ error is normalized by $\max_{\x} 4\pi\rho\approx 30.6796$. On the right are the associated convergence rates. The vertical line indicates the time at which the charge distribution stops moving.}  
        \label{fig:mc_stop_pics}
\end{figure*}

We have also performed another test with the same discretization and parameters but stopped the charge distribution after $t=\frac{40}{1024}$ and then run out to $t_{final}=\frac{100}{1024}$ to show that the solver recovers the electrostatic solution. Fig. \ref{fig:mc_stop_pics} shows the $E_x$ Richardson convergence rate estimate and the associated $\ell_{\infty}$ error as well as the absolute convergence rate and associated $\ell_{\infty}$ errors for $\nabla\cdot\E-4\pi\rho$ on the three grids in $\Omega_2$ as a function of time step and as expected our solution shows fourth order convergence. 

\begin{figure}[h!]
\centering
\includegraphics[scale=0.3]{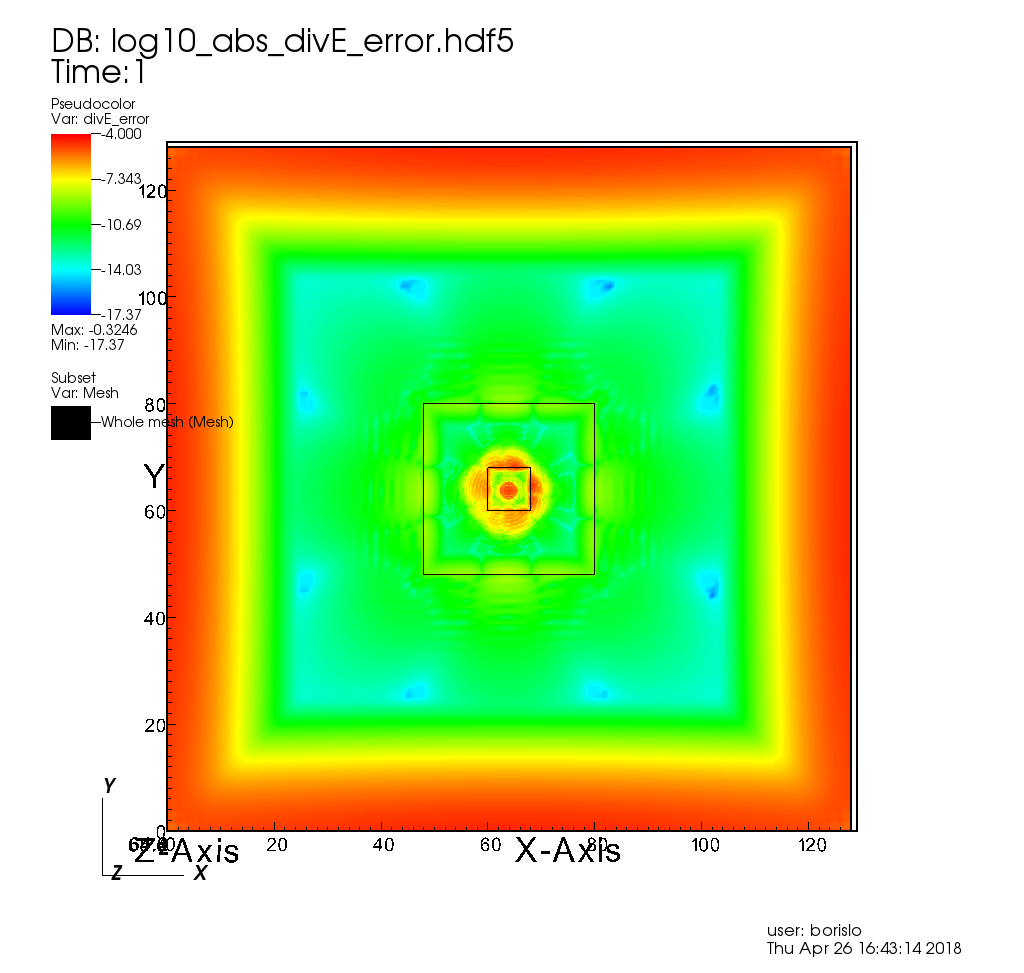}
\caption{$\log_{10}(|\nabla\cdot\E-4\pi\rho|/\max_{\x} 4\pi\rho)$ at $z=0.5$ for the stopped translating spherical charge distribution problem at $t=\frac{200}{2048}$ for $N=129$ showing that there are no reflected waves at the refinement boundaries.}
\label{fig:logabsplot}
\end{figure}

\subsubsection{Regridding test}
\begin{figure*}
        \centering
        \begin{subfigure}[t]{0.475\textwidth}
            \centering
            \includegraphics[width=\textwidth]{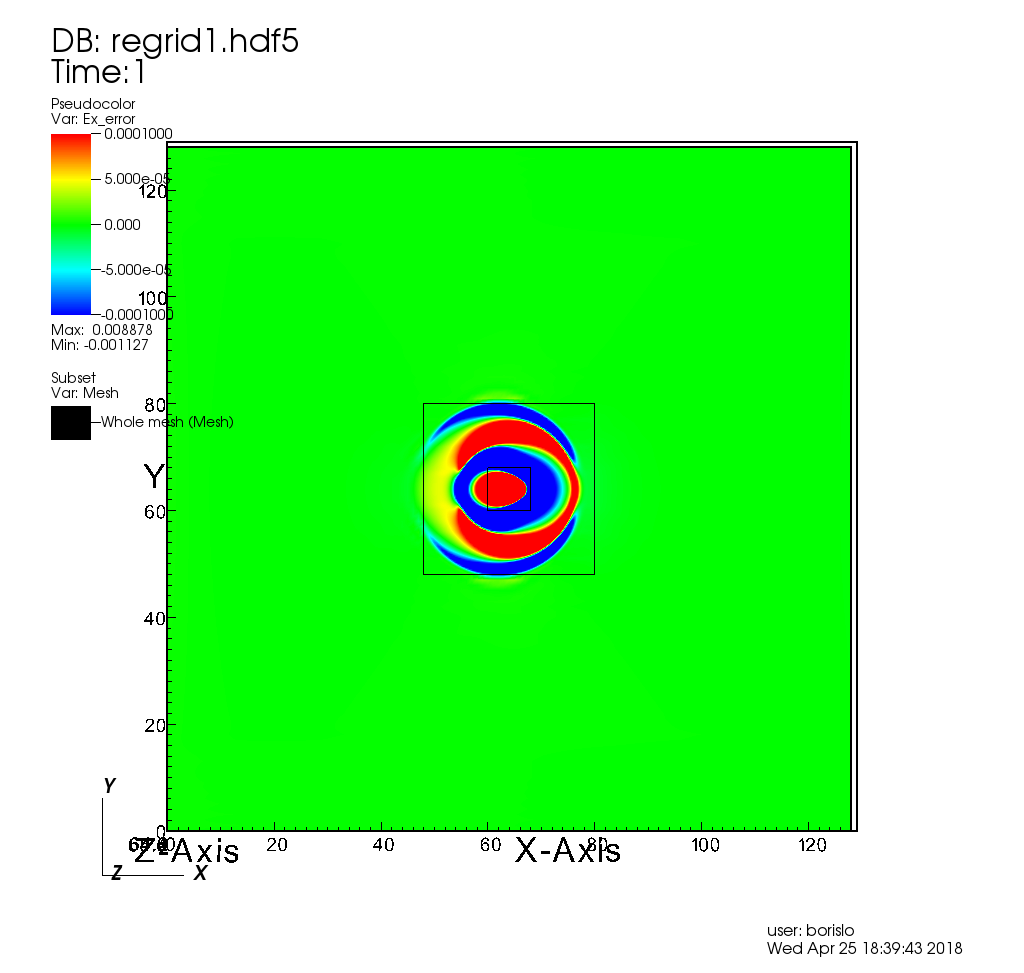}
            \caption[]%
            {{\small $t=\frac{256}{2048}$; charge distribution moving to the right, it has almost reached its rightmost position. $\Omega_2=\Omega_{2,b}$}}   
        \end{subfigure}
        \hfill
        \begin{subfigure}[t]{0.475\textwidth}  
            \centering 
            \includegraphics[width=\textwidth]{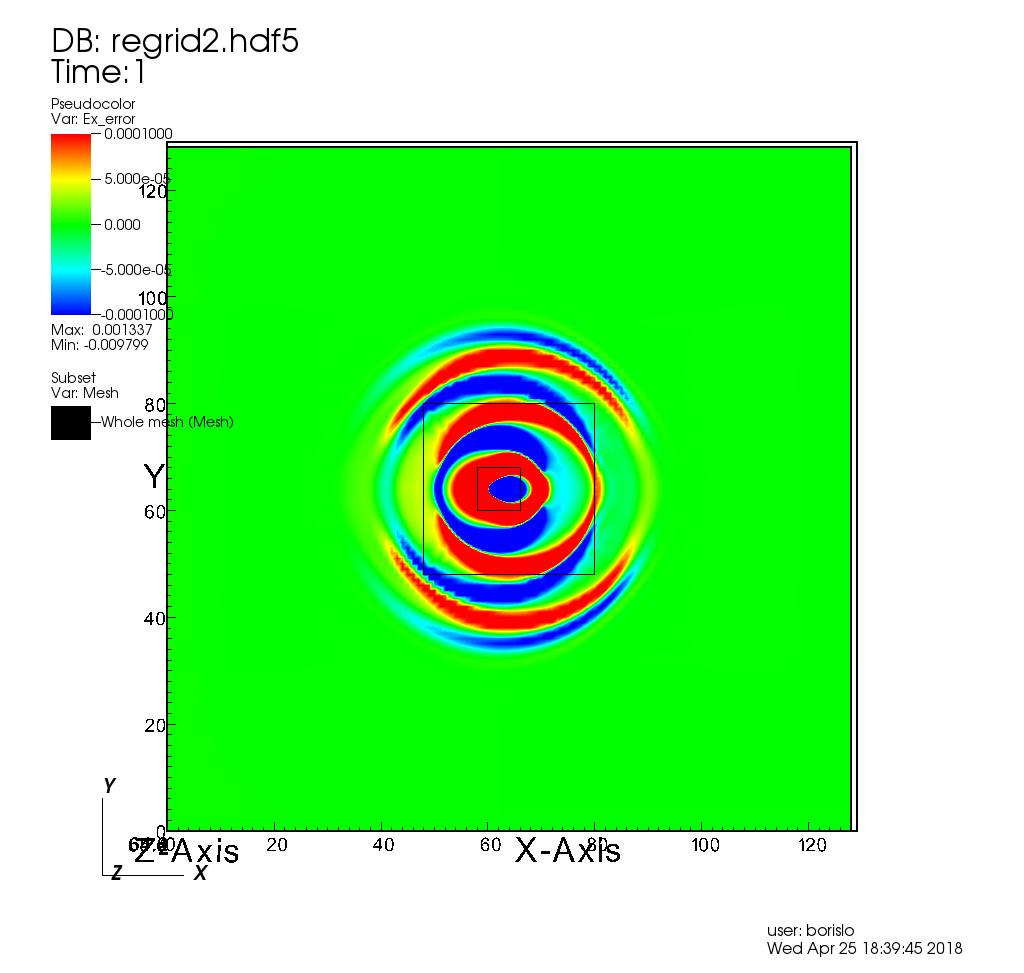}
            \caption[]%
            {{\small $t=\frac{480}{2048}$; charge distribution is at its leftmost position. $\Omega_2=\Omega_{2,a}$}}    
        \end{subfigure}
        \vskip\baselineskip
        \begin{subfigure}[t]{0.475\textwidth}   
            \centering 
            \includegraphics[width=\textwidth]{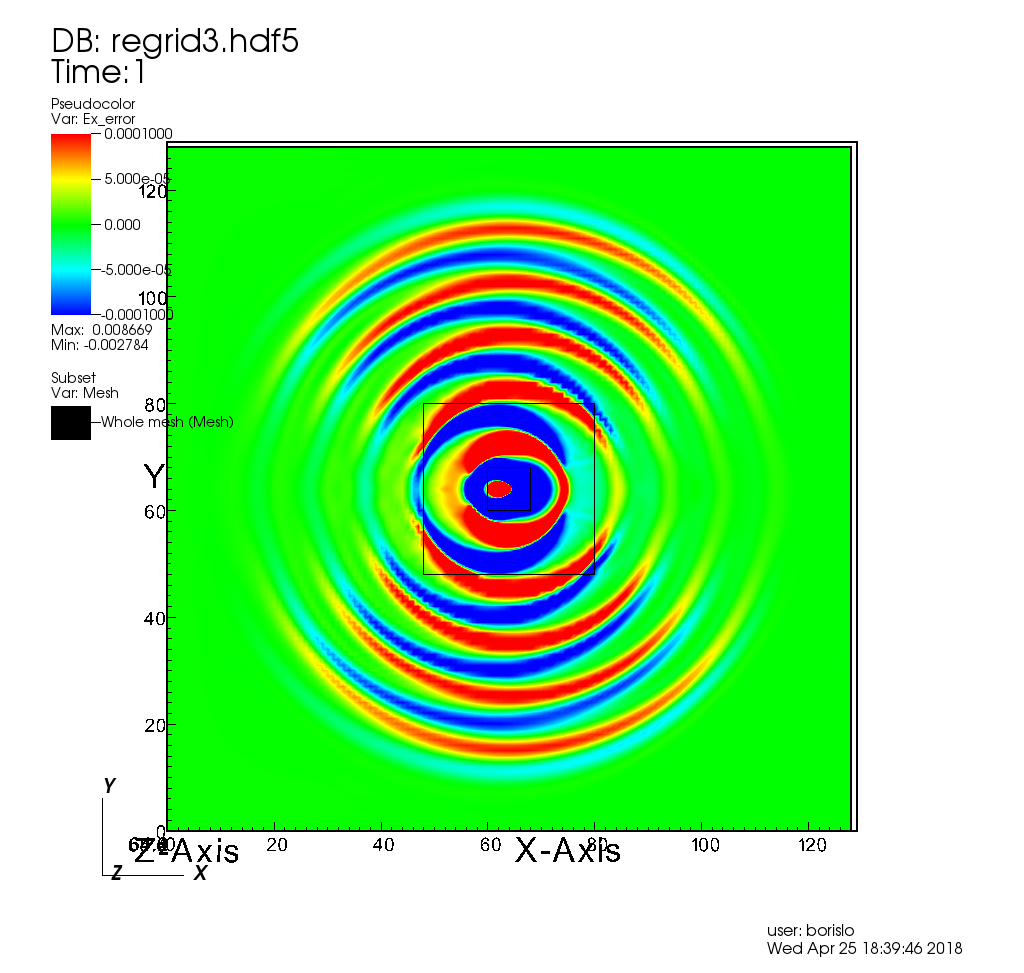}
            \caption[]%
            {{\small $t=\frac{864}{2048}$; charge distribution moving to the left. $\Omega_2=\Omega_{2,b}$}}    
        \end{subfigure}
        \quad
        \begin{subfigure}[t]{0.475\textwidth}   
            \centering 
            \includegraphics[width=\textwidth]{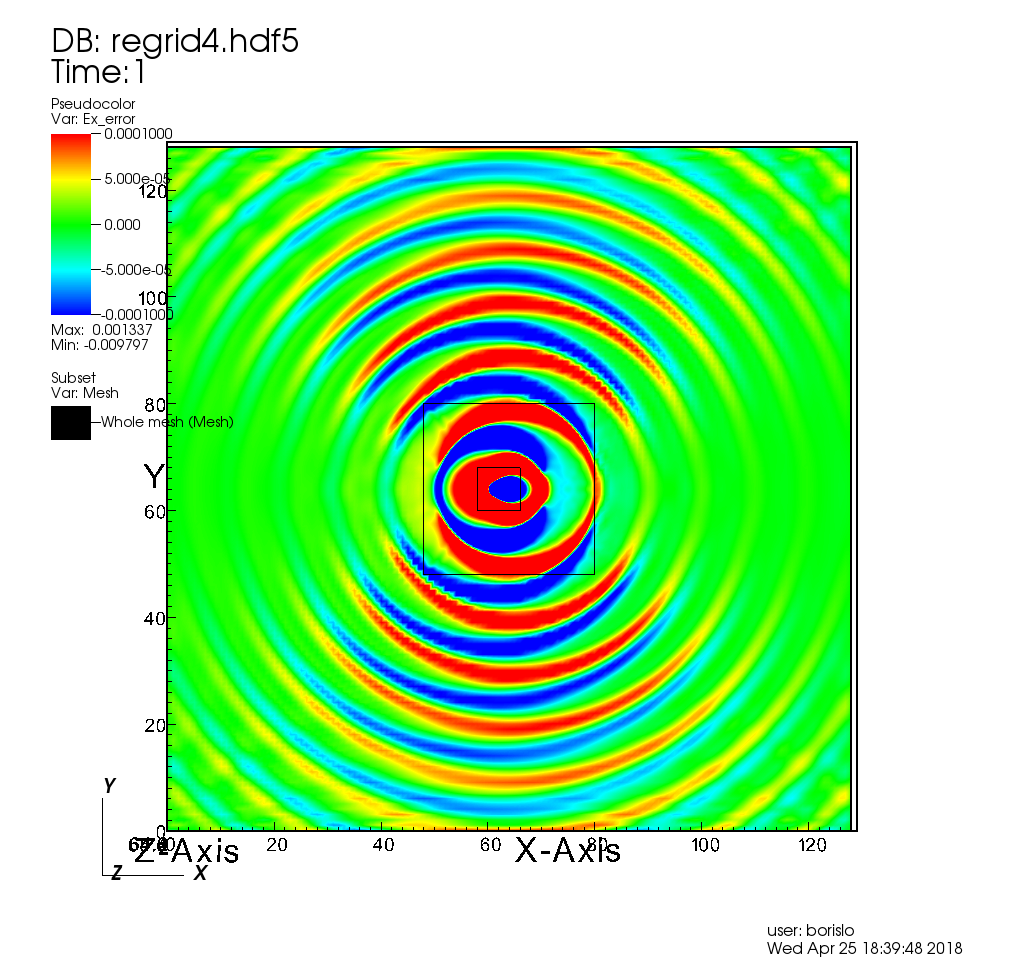}
            \caption[]%
            {{\small $t=\frac{1600}{2048}$; final time step. $\Omega_2=\Omega_{2,a}$}}    
        \end{subfigure}
        \caption[]
        {$E_x$ minus the instantaneous electrostatic solution, at $z = \frac{1}{2}$, for the spherical charge distribution problem with regridding for $N=129$.} 
        \label{fig:vx_regrid_pics}
\end{figure*}

\begin{figure*}
        \centering
        \begin{subfigure}[t]{0.475\textwidth}
            \centering
            \includegraphics[width=\textwidth]{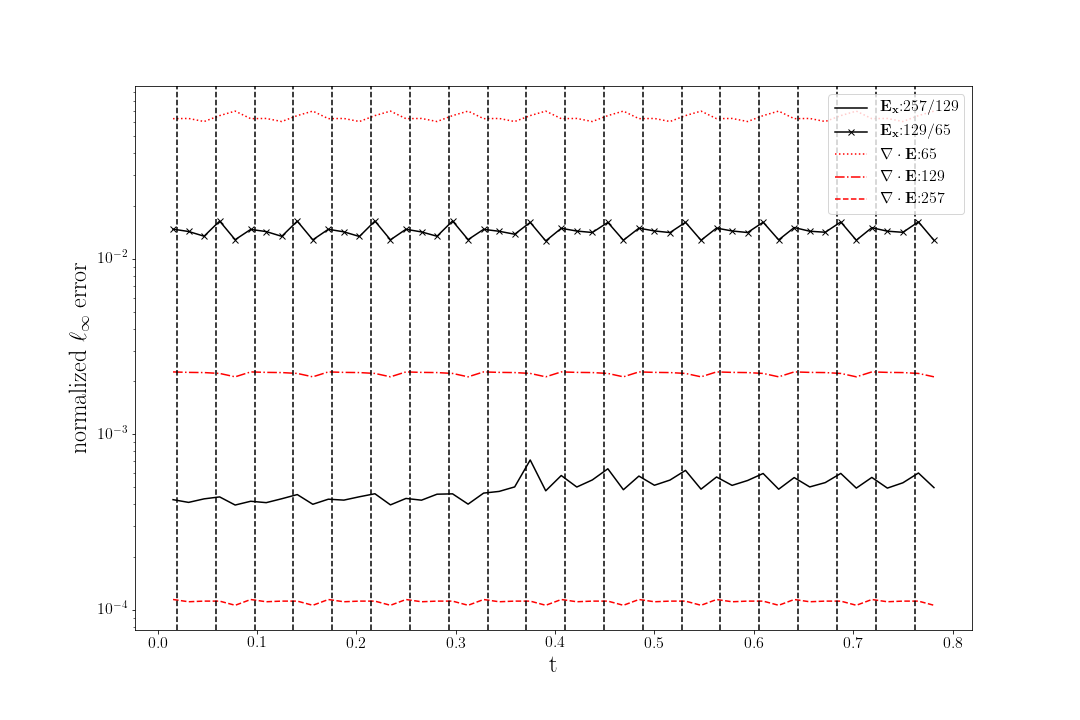} 
        \end{subfigure}
        \hfill
        \begin{subfigure}[t]{0.475\textwidth}  
            \centering 
            \includegraphics[width=\textwidth]{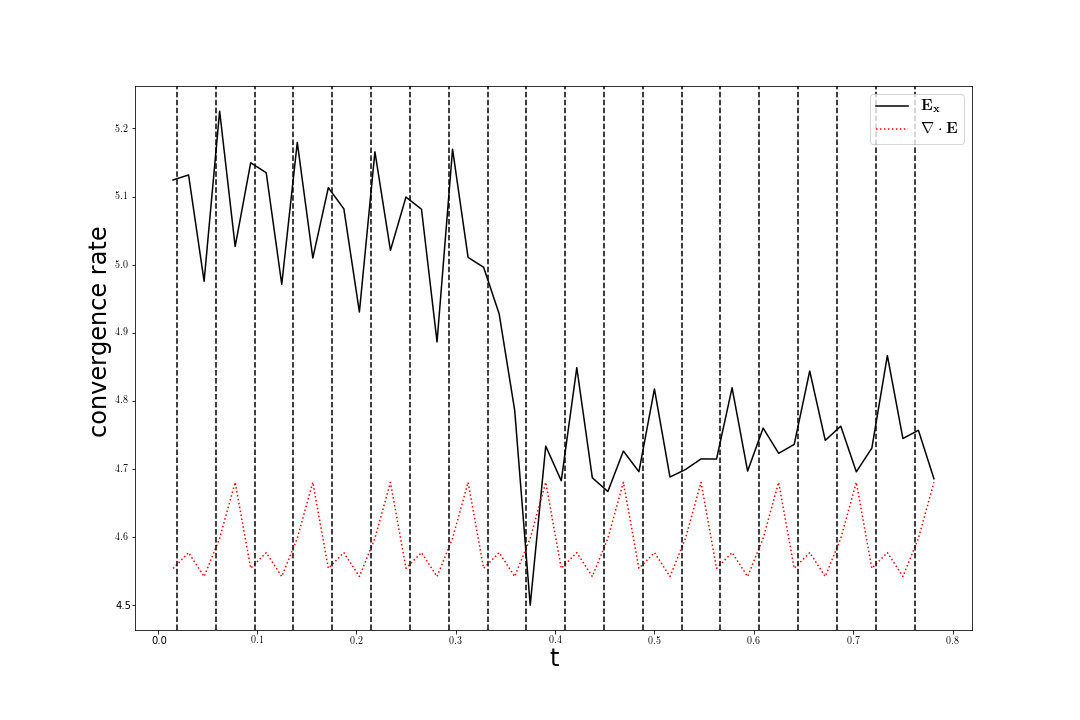}   
        \end{subfigure}
        \caption[]
         {$\ell_{\infty}$ error values and convergence results for $E_x$ and $\nabla\cdot\E-4\pi\rho$ for the regridding spherical charge distribution problem as a function of time in $\Omega_2$. On the left are the normalized $\ell_{\infty}$ errors for $E_x$ and $\nabla\cdot\E-4\pi\rho$. The errors for $E_x$ are obtained from the difference of sampled field values from the $N=257$ with $N=129$ and also from sampled $N=129$ with $N=65$ test case. The $E_x$ error is normalized by the max norm of the electrostatic solution ($\approx  0.0312658$) and $\nabla\cdot\E-4\pi\rho$ error is normalized by $\max_{\x} 4\pi\rho\approx 30.6796$. On the right are the associated convergence rates. The vertical lines are the times at which regridding occurs.}  
        \label{fig:regrid_pics}
\end{figure*}

We tested our regridding algorithm with the translating charge distribution with $\v=\nu d\pi \sin(2\nu t)\hat\x, a=\frac{1}{160}, d=\frac{1}{64}, \x_0 = \left(\frac{31}{64}, \frac{1}{2}, \frac{1}{2}\right), \nu = \frac{1024}{80}, t_{final}=\frac{800}{1024}$, and other parameters being the same. We kept $\Omega_1$ the same and fixed but regridded $\Omega_2$ starts with $\Omega_{2,a}$ and changes between $\Omega_{2,a}$ and $\Omega_{2,b}$ whenever the $x$ coordinate of the center of the charge distribution crosses $\frac{63}{128}$, where $\Omega_{2,a}$ is the rectangular prism defined by the corner points $(\frac{29}{64},\frac{17}{32},\frac{17}{32})$ and $(\frac{33}{64},\frac{17}{32},\frac{17}{32})$, and $\Omega_{2,b}=\left[\frac{15}{32},\frac{17}{32}\right]^3$; effectively $\Omega_2$ oscillates in the $x$ direction with amplitude $\frac{1}{64}$ in the direction of the charge motion. Fig \ref{fig:vx_regrid_pics} shows $E_x$ and the regridding domains for $N=129$.  Fig. \ref{fig:regrid_pics} shows the $E_x$ Richardson convergence rate estimate and the associated $\ell_{\infty}$ error as well as the absolute convergence rate and associated $\ell_{\infty}$ errors for $\nabla\cdot\E-4\pi\rho$ on the three grids in $\Omega_2$ as a function of time step and our solution shows fifth order convergence.  

\subsection{Divergence-Free Current Source}
\begin{figure*}
        \centering
        \begin{subfigure}[t]{0.475\textwidth}
            \centering
            \includegraphics[width=\textwidth]{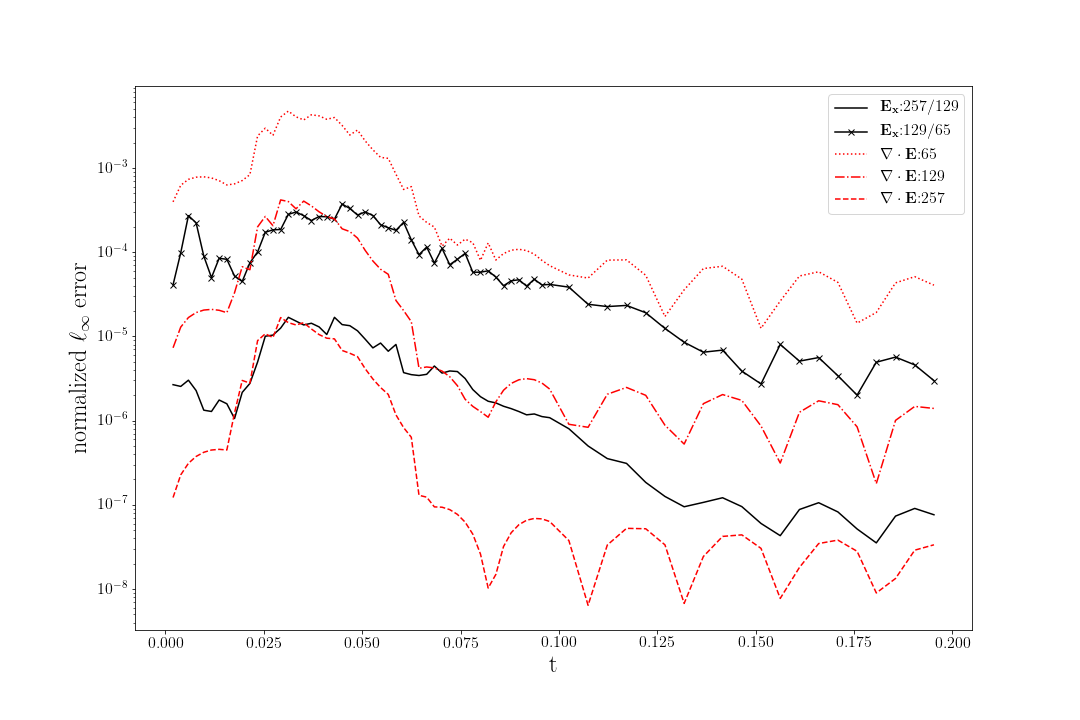}  
        \end{subfigure}
        \hfill
        \begin{subfigure}[t]{0.475\textwidth}  
            \centering 
            \includegraphics[width=\textwidth]{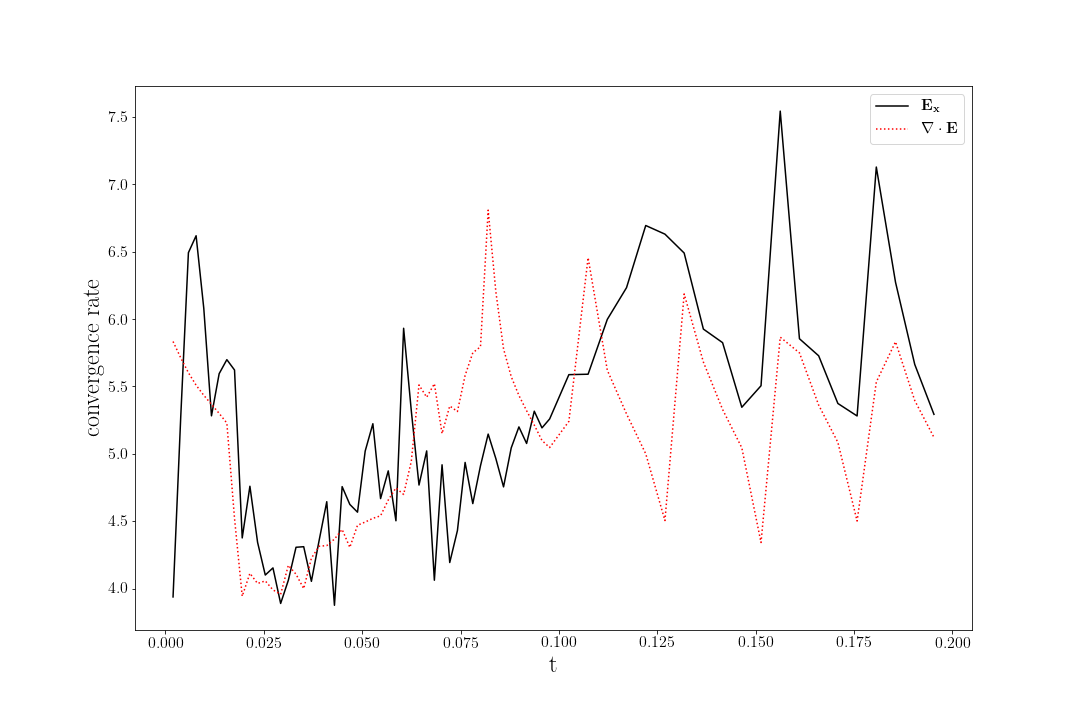} 
        \end{subfigure}
        \caption[]
          {$\ell_{\infty}$ error values and convergence results for $E_x$ and $\nabla\cdot\E-4\pi\rho$ for the divergence-free current problem as a function of time in $\Omega_2$. On the left are the normalized $\ell_{\infty}$ errors for $E_x$ and $\nabla\cdot\E-4\pi\rho$. The errors for $E_x$ are obtained from the difference of sampled field values from the $N=257$ with $N=129$ and also from sampled $N=129$ with $N=65$ test case. The $E_x$ error is normalized by $|\frac{4\pi}{\nu}\max_{r,z} J_x| \approx |10.2341\sin(2\pi\nu t)|$ and $\nabla\cdot\E$ is normalized by $|\frac{4\pi}{\nu a}\max_{r,z} J_x| \approx |545.8187\sin(2\pi\nu t)|$. On the right are the associated convergence rates.}  
        \label{fig:sven_pics}
\end{figure*}

We've also tested with a divergence-free current source of the form
\begin{align}
J_x (x,y,z,t) &= -100\frac{y-y_0}{r}\sin\left(\frac{\pi r}{2a}\right)\cos^{10}\left(\frac{\pi r}{2a}\right)\cos^{11}\left(\frac{\pi (z-z_0)}{d}\right)\sin(2\pi\nu t), \\
J_y (x,y,z,t) &= 100\frac{x-x_0}{r}\sin\left(\frac{\pi r}{2a}\right)\cos^{10}\left(\frac{\pi r}{2a}\right)\cos^{11}\left(\frac{\pi (z-z_0)}{d}\right)\sin(2\pi\nu t), \\
J_z (x,y,z,t) &= 0,
\end{align}
where $r=\sqrt{(x-x_0)^2 + (y-y_0)^2}$ with parameters: $a=\frac{3}{160}, d=\frac{13}{320}, x_0 = y_0 = z_0= 0.5, \nu = 20$, and using the same refinement levels and discretization, and $t_{final}$ as the fixed grids translating charge problem. Fig. \ref{fig:sven_pics} shows the $E_x$ Richardson convergence rate estimate and the associated $\ell_{\infty}$ error as well as the absolute convergence rate and associated $\ell_{\infty}$ errors for $\nabla\cdot\E$ on the three grids in $\Omega_2$ as a function of time step and as expected our solution shows fourth order convergence.

\section{Conclusion}\label{sec:conclusion}

We have presented a new version of our Green's function numerical method for Maxwell's equations. This new formulation results in a completely local propagator that does not require Helmholtz decomposition. In principle, the method can choose any CFL but at the cost of larger ghost regions. We have demonstrated a high order adaptive version of the solver in some test examples. In the future, we are interested in incorporating this method in EM PIC using Lawson's method where the fields and particles are evolved together with a Runge-Kutta scheme with an extra propagator step for the fields.

\section*{Acknowledgments}
This research is supported by the Office of Advanced Scientific Computing Research of the US Department of Energy under Contract Number DE-AC02-05CH11231. This research used resources of the National Energy Research Scientific Computing Center (NERSC), a Department of Energy Office of Science User Facility supported by the Office of Science of the U.S. Department of Energy under Contract Number DE-AC02-05CH11231. 

\begin{appendices}
\section{High Order B-Splines}\label{sec:appendix}

For completeness, we give the B-splines used in our implementation for the delta approximants and interpolants. Detailed discussions on creating high order B-splines are given in \cite{lo,chaniotis,walden,tornberg,schoenberg,monaghan}. $W_{q,p}$ denotes a $q$-th order accurate, $C^p$ B-spline.
\begin{enumerate}

\item 
\begin{align}
 W_{6,0}(x) &= \left\{
 \begin{array}{rl}
-\frac{|x|^5}{12 }+\frac{|x|^4}{4 }+\frac{5 |x|^3}{12 }-\frac{ 5 |x|^2}{4 }-\frac{ |x|}{3 }+ 1 & : |x|\in[0,1] \\[2pt]
\frac{|x|^5}{24 }-\frac{ 3 |x|^4}{8 }+\frac{25 |x|^3}{24 }-\frac{ 5 |x|^2}{8 }-\frac{ 13 |x|}{12 }+ 1 & : |x| \in[1,2] \\[2pt]
-\frac{|x|^5}{120 }+\frac{|x|^4}{8 }-\frac{ 17 |x|^3}{24 }+\frac{15 |x|^2}{8 }-\frac{ 137 |x|}{60 }+ 1 & : |x| \in[2,3]\\[2pt]
0 & : |x| > 3
\end{array}
 \right.
 \end{align}


\item 
\begin{align}
 W_{6,6}(x) &= \left\{
 \begin{array}{rl}
-\frac{665 |x|^9}{12048 }+\frac{ 665 |x|^8}{3012 }-\frac{ 2419 |x|^7}{12048 }-\frac{ 2437 |x|^6}{12048 } +\frac{ 2723 |x|^4}{3012 }-\frac{ 4543 |x|^2}{3012 }+\frac{ 19177}{21084} & :  |x|\in[0,1] \\[2pt]
\frac{133 |x|^9}{4016 }-\frac{ 399 |x|^8}{1004 }+\frac{ 39659 |x|^7}{20080 }-\frac{ 104409 |x|^6}{20080 } &\\[2pt]
+\frac{ 23443 |x|^5}{3012 }-\frac{ 14175 |x|^4}{2008 }+\frac{ 7553 |x|^3}{1506 }-\frac{ 32207 |x|^2}{10040 }+\frac{ 2933 |x|}{15060 }+\frac{ 13081}{14056} & :  |x|\in[1,2] \\[2pt]
-\frac{133 |x|^9}{12048 }+\frac{ 665 |x|^8}{3012 }-\frac{ 114139 |x|^7}{60240 }+\frac{ 109283 |x|^6}{12048 } &\\[2pt]
-\frac{ 79303 |x|^5}{3012 }+\frac{ 283423 |x|^4}{6024 }-\frac{ 75215 |x|^3}{1506 }+\frac{ 170023 |x|^2}{6024 }-\frac{ 90923 |x|}{15060 }-\frac{ 17653}{42168} & :  |x|\in[2,3] \\[2pt]
\frac{19 |x|^9}{12048 }-\frac{ 133 |x|^8}{3012 }+\frac{ 225859 |x|^7}{421680 }-\frac{ 221003 |x|^6}{60240 } &\\[2pt]
+\frac{ 23299 |x|^5}{1506 }-\frac{ 30793 |x|^4}{753 }+\frac{ 49184 |x|^3}{753 }-\frac{ 208208 |x|^2}{3765 }+\frac{ 53632 |x|}{3765 }+\frac{ 32512}{5271} & :  |x|\in[3,4] \\[2pt]
0 & : |x| > 4
\end{array}
 \right.
 \end{align}
\end{enumerate}
\end{appendices}

\bibliographystyle{ieeetr}
\bibliography{reference}

\end{document}